\documentclass[11pt]{amsart}

\theoremstyle{plain}
\newtheorem{thm}{Theorem}
\newtheorem{theorem}[thm]{Theorem}

\theoremstyle{definition}
\newtheorem{defn}[thm]{Definition}

\theoremstyle{remark}

\newtheorem*{remark*}{Remark}

\begin{document}

\author{Andrey A. Yakovlev}
\address{Department of Mathematics, Ufa State Aviation Technical University,
12 K. Marx str., 450000 Ufa, Russia}
\email{yakovlevandrey@yandex.ru}

\title{Adiabatic limits on Riemannian Sol-manifolds}
\thanks{Supported by the Russian Foundation of Basic Research
(grant no. 06-01-00208)}

\begin{abstract}
We obtain an asymptotic formula for the spectrum distribution
function of the Laplace operator on a compact Riemannian
Sol-manifold in the adiabatic limit determined by a one-dimensional
foliation defined by the orbits of a left-invariant flow.
\end{abstract}

\maketitle

The paper is devoted to investigation of adiabatic limits on
Riemannian Sol-manifolds. We understand adiabatic limits in the
sense, which was introduced by Witten in~\cite{Witten}. More
precisely, let $(M,{\mathcal F})$ be a closed foliated manifold
equipped with a Riemannian metric $g$. Thus, the tangent bundle $TM$
of $M$ is represented as a direct sum
\[
TM=F\oplus H,
\]
where $F=T{\mathcal F}$ is the tangent bundle of $\mathcal F$ and
$H=F^{\bot}$ the orthogonal complement of $F$. Let $g_F$ and $g_H$
denote the restriction of the metric $g$ to $F$ and $H$,
respectively. Therefore, $g=g_F+g_H$. Define a one-parameter family
of Riemannian metrics on $M$ by the formula
\begin{equation}\label{e:gh}
g_{\varepsilon}=g_{F} + {\varepsilon}^{-2}g_{H}, \quad \varepsilon
> 0.
\end{equation}
Investigation of various properties of the family of Riemannian
manifolds $(M, g_\varepsilon)$ as $\varepsilon\to 0$ will be called
by passage to adiabatic limit.

Recall \cite{Thurston} that the group $Sol$ is the solvable Lie
subgroup of the Lie group $\operatorname{GL}(3,\mathbb{R})$, which
consists of all matrices of the form:
\[
\gamma(u,v,w)=\begin{pmatrix}
 e^w & 0      &u\\
 0   & e^{-w} &v\\
 0   & 0      &1
\end{pmatrix},\quad (u,v,w)\in \mathbb R^3.
\]
The Lie algebra $\mathfrak{sol}$ of $Sol$ is the Lie subalgebra of
the Lie algebra $gl(3,\mathbb{R})$, which consists of all matrices
of the form
\[
X(u,v,w)=\begin{pmatrix}
 w & 0      &u\\
 0   & -w &v\\
 0   & 0      &0
\end{pmatrix},\quad (u,v,w)\in \mathbb R^3.
\]

Let $A\in\operatorname{SL}(2,\mathbb{Z})$ and $|\operatorname {tr}
A|>2$. Denote by $\lambda$ and $\lambda^{-1}$ the eigenvalues of $A$
and assume that $\lambda>1$. Define a vectors
$(c^1_1,c^2_1),(c^1_2,c^2_2)$ by the equation
\[
A=\begin{pmatrix} a_{11} & a_{12}\\ a_{21} & a_{22}
\end{pmatrix}=\begin{pmatrix} c^1_1 & c^1_2\\ c^2_1 & c^2_2
\end{pmatrix}^{-1}\begin{pmatrix} \lambda & 0\\ 0 & \lambda^{-1}
\end{pmatrix}\begin{pmatrix} c^1_1 & c^1_2\\ c^2_1 & c^2_2
\end{pmatrix}.
\]

\begin{defn} A Riemannian Sol-manifold is a oompact manifold
 $M^3_A=G_A\backslash Sol$ equipped with a Riemannian metric $g$, where:
\begin{itemize}
\item $G_A$ is the uniform discrete subgroup of the Lie group
$Sol$, which consists of all $\gamma(u,v,w)\in Sol$ such that
  \[
(u,v)\in\Gamma:=\{ k(c^1_1,c^2_1)+l(c^1_2,c^2_2), \quad k,l\in
\mathbb Z\},
  \]
  \[
\quad w=m\ln\lambda, \quad m\in \mathbb Z,
  \]
\item $g$ is a Riemannian metric on $M^3_A$ whose lift on
$Sol$ is invariant under left translations by elements of $Sol$
(such metrics will be called locally left-invariant).
\end{itemize}
\end{defn}

A locally left-invariant metric $g$ is uniquely determined by its
value at the identity $\gamma(0,0,0)$ of $Sol$, and, therefore, is
given by a symmetric positive definite $3\times 3$-matrix.

Let $\alpha  \in \mathbb{R}$. Consider the left-invariant vector
field on $Sol$ associated with $X(1,\alpha,0)\in \mathfrak{sol}$.
The orbits of the corresponding vector field on $M^3_A$ define a
one-dimensional foliation $\mathcal F$. The leaf of $\mathcal F$
through $G_A \gamma(u,v,w)\in M^3_{A}$ is given by
\[
L_{G_A \gamma(u,v,w)} =\{ G_A\gamma(u+e^w t,v+\alpha e^{-w} t,w)\in
 M^3_{A}: t\in {\mathbb R}\}.
\]

Suppose that a locally left-invariant metric $g$ correspond to the
identity matrix. Consider the adiabatic limit associated with the
Riemannian Sol-manifold $(M^3_{A},g)$ and the foliation $\mathcal
F$. Denote by $\Delta_\varepsilon$ the Laplace-Beltrami operator on
$ M^3_{A}$ associated with the metric $g_{\varepsilon}$ given by
(\ref{e:gh}). For any $\varepsilon>0$ the spectrum of
$\Delta_\varepsilon$ consists of eigenvalues of finite multiplicity:
\[
0=\lambda_0(\varepsilon)< \lambda_1(\varepsilon)\leq \ldots,
\lambda_j(\varepsilon)\to +\infty\ \text{при}\ j\to\infty.
\]

The main result of the paper is a computation of the asymptotics of
the spectrum distribution function
\begin{equation*}
N_\varepsilon(t)=\sharp \{i:\lambda_i(\varepsilon)\leq t\}
\end{equation*}
of the operator $\Delta_\varepsilon$ in the adiabatic limit, that
is, when $t\in\mathbb{R}$ is fixed and $\varepsilon\rightarrow 0$.

\begin{theorem}
For any $t>0$, the following asymptotic formulae hold:
\medskip\par
1. For $\alpha\neq 0$
\begin{equation*}
N_\varepsilon(t)=\frac{1}{4\pi^2}t^{\frac{3}{2}}\varepsilon^{-2}+o(\varepsilon^{-2}),\quad
\varepsilon\to0.
\end{equation*}
\medskip\par
2. For $\alpha= 0$
\begin{equation*}
N_\varepsilon(t)=\frac{1}{6\pi^2}t^{\frac{3}{2}}\varepsilon^{-2}+o(\varepsilon^{-2}),\quad
\varepsilon\to0.
\end{equation*}
\end{theorem}

The asymptotic behavior of the spectrum distribution function for
the Laplace operator in the adiabatic limit was studied earlier in
\cite{adiab} for Riemannian foliations and in \cite{Ya_MS} for
one-dimensional foliations on Riemannian Heisenberg manifolds (see
also \cite{Bedlewo}). In all cases, the function $N_\varepsilon (t)$
has order $\varepsilon^{-q}$, where $q$ is the codimension of the
foliation (in our case $q=2$), but the coefficients of
$\varepsilon^{-q}$ are different in each case. Observe also that, in
each case, the asymptotic formula for $N_\varepsilon(t)$ in the
adiabatic limit is different from the classical Weyl formula, which
describes asymptotic behavior of $N_\varepsilon(t)$ as $t\to \infty$
(cf. \cite{Bolsinov}).

For $\alpha\neq 0$, the proof of the theorem uses the calculation of
the spectrum of the Laplace operator on a Riemannian Sol-manifold
given in \cite{Bolsinov}, which continues the investigation of the
geodesic flow on a Riemannian Sol-manifold started in
\cite{TaimanovBolsinov1} and \cite{TaimanovBolsinov2}, and
semiclassical spectral asymptotics \cite{Helffer} for the modified
Mathieu operator
\[
H_\varepsilon=-\varepsilon^2\frac{d^2}{dx^2}+a\operatorname{ch}(2\mu
x), \quad x\in {\mathbb R}.
\]

In the case $\alpha= 0$, the foliation is Riemannian, and the metric
is bundle-like, and, therefore, we can use the asymptotic formula
obtained in \cite{adiab}.

The author is grateful to Yu.A. Kordyukov for posing this problem
and attention to his work and to I.A. Taimanov for useful remarks.


\end{document}